\documentclass[11pt]{article}

\usepackage{amsmath}
\usepackage{mathrsfs}
\usepackage{graphicx}
\usepackage{amssymb}
\usepackage{a4}
\usepackage{url}
\usepackage{natbib}

\def\dd{{\rm {d}}}

\def\mp{\partial}
\def\tr{\mbox{\rm tr}}
\def\Ib{{\bf I}}
\def\mt{\theta}

\def\me{\epsilon}

\def\ra{\rightarrow}

\def\TT{^\top}

\def\ml{\lambda}
\def\ma{\alpha}
\def\mo{\omega}
\def\mob{\mathbf{\omega}}

\def\mg{\gamma}
\def\Ab{\mathbf{A}}
\def\Bb{\mathbf{B}}
\def\Hb{\mathbf{H}}
\def\mLb{\mathbf{\Lambda}}
\def\Mb{\mathbf{M}}
\def\Sb{\mathbf{S}}
\def\sb{\mathbf{s}}
\def\xb{\mathbf{x}}
\def\vb{\mathbf{v}}

\def\zb{\mathbf{z}}

\def\0b{\mbox{\rm 0}}

\def\SE{{\mathscr E}}

\def\SM{{\mathscr M}}

\def\SX{{\mathscr X}}

\newtheorem{rmk}{Remark}
\newtheorem{thm}{Theorem}
\newtheorem{lemma}{Lemma}

\newcommand{\vsp}{\vspace{0.3cm}}

\begin{document}

\title{A delimitation of the support of optimal designs \\for Kiefer's $\phi_p$-class of criteria}

\author{{L. {\sc Pronzato}} \\
\mbox{}\\
Laboratoire I3S, CNRS/Universit\'e de Nice-Sophia Antipolis \\
B\^at.\ Euclide, Les Algorithmes, 2000 route des Lucioles, BP 121\\
06903 Sophia Antipolis cedex, France \\
%\email{
{\tt pronzato@i3s.unice.fr}
}

\maketitle %[P]

\begin{abstract}
The paper extends the result of Harman and Pronzato [Stat. \& Prob. Lett., 77:90--94, 2007], which corresponds to $p=0$, to all strictly concave criteria in Kiefer's $\phi_p$-class. We show that, for any given design measure $\xi$, any support point $\xb_*$ of a $\phi_p$-optimal design
is such that the directional derivative of $\phi_p$ at $\xi$ in the direction of the delta measure at $\xb_*$ is larger than some bound $h_p[\xi]$ which is easily computed: it requires the determination of the unique root of a simple univariate equation (polynomial when $p$ is integer) in a given interval. The construction can be used to accelerate algorithms for $\phi_p$-optimal design and is illustrated on an example with $A$-optimal design.
\end{abstract}

{\bf keywords} Approximate design; optimum design; support points; design algorithm


%______________________________________________________________________________________
\section{Introduction and motivation}\label{S:Intro}

For $\SX$ a compact subset of $\mathbb{R}^m$, denote by $\Xi$ the set of design measures (i.e., probability measures) on $\SX$ and by $\Mb(\xi)$ the information matrix
$$
    \Mb(\xi)=\int_\SX \xb\xb\TT\,\xi(\dd \xb) \,.
$$
We suppose that there exists a nonsingular design on $\SX$ (i.e., there exists a $\xi\in\Xi$ such that $\Mb(\xi)$ is nonsingular) and we denote by $\Xi^+$ the set of such designs.
We consider an optimal design problem on $\SX$ defined by the maximization of a design criterion $\phi(\xi)=\Phi[\Mb(\xi)]$ with respect to $\xi\in\Xi$.
One may refer to \citet[Chap.~5]{Pukelsheim93} for a presentation of desirable properties that make a criterion $\Phi(\cdot)$ appropriate to measure the information provided by $\xi$. Here we shall focus our attention on design criteria that correspond to the $\phi_p$-class considered by \citet{Kiefer74}. More precisely, we consider the positively homogeneous form of such criteria and, for any $\Mb\in\mathbb{M}$, the set of symmetric non-negative definite $m\times m$ matrices, we denote
\begin{equation}\label{Phi_p^+}
    \Phi_p^+(\Mb)= \left[\frac1m\, \tr (\Mb^{-p}) \right]^{-1/p} \,, %\ p\in(-1,\infty)
\end{equation}
with the continuous extension $\Phi_p^+(\Mb)=0$ when $\Mb$ is singular and $p \geq 0$.
A design measure $\xi_p^*$ that maximizes $\phi_p(\xi)=\Phi_p^+[\Mb(\xi)]$ will be said $\phi_p$-optimal.
Note that when $p\neq 0$ the maximization of $\Phi_p^+(\Mb)$ is equivalent to the minimization of $\left[\tr (\Mb^{-p}) \right]^{1/p}$, and thus to the minimization of $\tr (\Mb^{-p})$ when $p$ is positive. A classical example is $A$-optimal design, which corresponds to $p=1$. Taking the limit of $\Phi_p^+(\cdot)$ when $p$ tends to zero, we obtain
$
\Phi_0^+(\Mb)= [\det(\Mb)]^{1/m} \,,
$
which corresponds to $D$-optimal design. The limit when $p$ tends to infinity gives $\Phi_\infty(\Mb)=\ml_{\min}(\Mb)$, the minimum eigenvalue of $\Mb$, and corresponds to $E$-optimal design. Some basic properties of $\phi_p$-optimal designs are briefly recalled in Sect.~\ref{S:properties Phi_p}.
%A key property among those is that a finitely-supported optimal design always exists.

Classical algorithms for optimal design usually apply to situations where $\SX$ is a finite set. The performance of the algorithm (in particular, its execution time for a given required precision on $\phi(\cdot)$) then heavily depends on the number $k$ of elements in $\SX$. The case of $D$-optimal design has retained much attention, see, for instance, \citet{AhipasaogluST2008}, \citet{ToddY2007}, \citet{Yu2010} and \citet{Yu2010-SC}. \citet{HPa06_SPL} show how any nonsingular design on $\SX$ yields a simple inequality that must be satisfied by the support points of a $D$-optimal design $\xi_0^*$. Whatever the iterative method used for the construction of $\xi_0^*$, this delimitation of the support of $\xi_0^*$ permits to reduce the cardinality of $\SX$ along the iterations, with the inequality becoming more stringent when approaching the optimum, hence producing a significant acceleration of the algorithm. Put in other words, the delimitation of the support of an optimal design facilitates the optimization by focussing the search on the useful part of the design space $\SX$. The objective of the paper is to extend the results in \citet{HPa06_SPL} to the $\phi_p$-class (\ref{Phi_p^+}) of design criteria. The condition obtained does not tell what the optimum support is, but indicates where it cannot be.

The paper is organized as follows. Section~\ref{S:properties Phi_p} recalls the main properties of $\phi_p$-optimal design that are useful for the rest of the paper. The main result is derived in Sect.~\ref{S:Phi_p} and illustrative examples are given in Sect.~\ref{S:Example}. Finally, Sect.~\ref{S:Extensions} concludes and indicates some possible extensions. The technical parts of the proofs are given in appendix.

%______________________________________________________________________________________
\section{Some basic properties of $\phi_p$-optimal designs}\label{S:properties Phi_p}

The criteria $\Phi_p^+(\cdot)$ defined by (\ref{Phi_p^+}) satisfy $\Phi_p^+(\Ib_m)=1$ for $\Ib_m$ the $m$-dimensional identity matrix and $\Phi_p^+(a\Mb)=a\,\Phi_p^+(\Mb)$ for any $a>0$ and any $\Mb\in\mathbb{M}$.
Note that, from Caratheodory's theorem, a finitely-supported optimal design always exists, with $m(m+1)/2$ support points at most. We also have the following properties.

\begin{lemma}\label{L:properties} For any $p\in(-1,\infty)$, the criterion $\Phi_p^+(\cdot)$ satisfies the following:
\begin{itemize}
  \item[(i)] $\Phi_p^+(\cdot)$ is strictly concave on the set $\mathbb{M}^+$ of symmetric positive definite $m\times m$ matrices; it is strictly isotonic (it preserves L{\"o}wner ordering) on $\mathbb{M}$ for $p\in(-1,0)$; that is, $\Phi_p^+(\Mb_2)>\Phi_p^+(\Mb_1)$ for all $\Mb_1$ and $\Mb_2$ in $\mathbb{M}$ such that $\Mb_2-\Mb_1\in\mathbb{M}$ and $\Mb_2 \neq \Mb_1$;
  it is strictly isotonic on $\mathbb{M}^+$ for $p\in[0,\infty)$.
  \item[(ii)] Any $\phi_p$-optimal design $\xi_p^*$ is nonsingular.
  \item[(iii)] The optimal matrix $\Mb_*=\Mb_*[p]$ is unique.
\end{itemize}
\end{lemma}

Part ($i$) is proved in \citet[Chap.~6]{Pukelsheim93}.
For $p\geq 0$, ($ii$) follows from the observation that $\Phi_p^+(\Mb)=0$ when $\Mb$ is singular while there exists a nonsingular $\Mb(\xi)$ with $\Phi_p^+[\Mb(\xi)]>0$;
for $p\in(-1,0)$, the statement is proved in \citep[Sect.~7.13]{Pukelsheim93} through the use of polar information functions.
Part ($iii$) is a direct consequence of ($i$) and ($ii$): since an optimal design matrix $\Mb_*$ is nonsingular, the strict concavity of $\Phi_p^+(\cdot)$ at $\Mb_*$ implies that $\Mb_*$ is unique. Note that this does not imply that the optimal design measure $\xi_p^*$ maximizing $\phi_p(\xi)$ is unique.

\vsp
We shall only consider values of $p$ in $(-1,\infty)$ and, from Lemma~\ref{L:properties}-($ii$), we can thus restrict our attention to matrices $\Mb$ in $\mathbb{M}^+$. $\Phi_p^+(\cdot)$ is differentiable at any $\Mb\in\mathbb{M}^+$, with gradient
$$
\nabla\Phi_p^+(\Mb) = \frac1m\, [\Phi_p^+(\Mb)]^{p+1}\, \Mb^{-(p+1)} = \frac{\Phi_p^+(\Mb)}{\tr (\Mb^{-p})}\, \Mb^{-(p+1)} \,.
$$
The directional derivative
$F_{\phi_p}(\xi;\nu)=\lim_{\ma\ra 0^+} (1/\ma) \{\phi_p[(1-\ma)\xi+\ma\nu]-\phi_p(\xi)\}$
is well defined and finite for any $\xi\in\Xi^+$ and any $\nu\in\Xi$, with
$$
F_{\phi_p}(\xi;\nu) = \tr\{[\Mb(\nu)-\Mb(\xi)]\nabla\Phi_p^+[\Mb(\xi)]\}
= \phi_p(\xi)\, \left\{ \frac{\int_\SX \xb\TT \Mb^{-(p+1)}(\xi)\xb \,\nu(\dd\xb)}{\tr [\Mb^{-p}(\xi)]} - 1 \right\} \,.
$$
We shall denote by $F_{\phi_p}(\xi,\xb)=F_{\phi_p}(\xi;\delta_\xb)$ the directional derivative of $\phi_p(\cdot)$ at $\xi$ in the direction of the delta measure at $\xb$,
\begin{equation}\label{F(xi,x)}
    F_{\phi_p}(\xi,\xb) = \phi_p(\xi)\, \left\{ \frac{\xb\TT \Mb^{-(p+1)}(\xi)\xb}{\tr [\Mb^{-p}(\xi)]} - 1 \right\} \,.
\end{equation}

The following theorem, which relies on the concavity and differentiability of $\Phi_p^+(\cdot)$, is a classical result in optimal design theory, see, e.g., \citet{Kiefer74} and \citet[Chap.~7]{Pukelsheim93}.

\begin{thm}[Equivalence Theorem]\label{Thm:EQ}
For any $p\in(-1,\infty)$, the following statements are equivalent:
\begin{itemize}
    \item[(i)] $\xi_p^*$ is $\phi_p$-optimal.
    \item[(ii)] $\xb\TT \Mb^{-(p+1)}(\xi_p^*)\xb \leq \tr [\Mb^{-p}(\xi_p^*)]$ for all $\xb\in\SX$.
    \item[(iii)] $\xi_p^*$ minimizes $\max_{\xb\in\SX} F_{\phi_p}(\xi,\xb)$ with respect to $\xi\in\Xi^+$.
\end{itemize}
Moreover, the inequality of ($ii$) holds with equality for every support point $\xb=\xb_*$ of $\xi_p^*$.
\end{thm}

%______________________________________________________________________________________
%\section{Support delimitation for Kiefer's $\Phi_p$ class of optimality criteria}\label{S:Phi_p}
\section{A necessary condition for support points of $\phi_p$-optimal designs}\label{S:Phi_p}

%--------------------------------------------------------------------------------------
\subsection{A lower bound on $\xb\TT\Mb^{-(p+1)}\xb$ for the support points of an optimal design}

Take any $p\in(-1,\infty)$ and any $\xi\in\Xi^+$. We shall omit the dependence in $\xi$ when there is no ambiguity and simply write $\Mb=\Mb(\xi)$, $\phi_p=\phi_p(\xi)$. We shall also denote
$$
t=t(\xi,p)=\tr [\Mb^{-p}]\,, \ \ t_*=t_*(p)=\tr(\Mb_*^{-p})\,,
$$
with $\Mb_*$ the optimal matrix satisfying $\phi_p^*=\Phi_p^+(\Mb_*)=\max_{\nu\in\Xi} \Phi_p^+[\Mb(\nu)]$. Define
\begin{equation}\label{epsilon}
    \me=\me(\xi,p)=\max_{\xb\in\SX} \{\xb\TT \Mb^{-(p+1)}\xb\} - t \,.
\end{equation}
The concavity of $\Phi_p^+(\cdot)$ implies that
$
\phi_p \leq \phi_p^* \leq \phi_p+ F_{\phi_p}(\xi;\xi_p^*) \leq \phi_p(\xi)+ \max_{x\in\SX} F_{\phi_p}(\xi,\xb)\,,
$
with $\xi_p^*$ denoting a $\phi_p$-optimal design measure; that is,
\begin{equation}\label{ineq-phi-phi_*-me}
    \phi_p \leq \phi_p^* \leq \phi_p\,(1+\me/t)\,,
\end{equation}
see (\ref{F(xi,x)}).

Since $\xb\TT \Mb^{-(p+1)}\xb \leq t+\me$ for all $\xb\in\SX$, see (\ref{epsilon}), we have
\begin{equation}\label{ineq1-t+me}
\tr[\Mb_*\Mb^{-(p+1)}] \leq t+\me \,.
\end{equation}
On the other hand, the optimality of $\xi_p^*$ implies (see Th.~\ref{Thm:EQ}-($ii$))
\begin{equation}\label{ineq2-t_*}
\tr[\Mb\Mb_*^{-(p+1)}] \leq t_* \,.
\end{equation}
Moreover, any support point $\xb_*$ of $\xi_p^*$ satisfies $\xb_*\TT \Mb_*^{-(p+1)}\xb_* = t_*$. We use a construction similar to that in \citet{HPa06_SPL} and define $\Hb=\Hb(\xi,p)=\Mb^{-(p+1)/2}\Mb_*^{p+1}\Mb^{-(p+1)/2}$. Then we can write
%\begin{eqnarray*}
%\xb_*\TT \Mb^{-(p+1)}\xb_* &=& \xb_*\TT \Mb^{-(p+1)/2} \Hb^{-1/2} \Hb \Hb^{-1/2} \Mb^{-(p+1)/2} \xb_* \\
%&& \geq  \ml_1\, \xb_*\TT \Mb_*^{-(p+1)} \xb_* = \ml_1\, t_* \,.
%\end{eqnarray*}
$$
\xb_*\TT \Mb^{-(p+1)}\xb_* = \xb_*\TT \Mb^{-(p+1)/2} \Hb^{-1/2} \Hb \Hb^{-1/2} \Mb^{-(p+1)/2} \xb_* \geq  \ml_1\, \xb_*\TT \Mb_*^{-(p+1)} \xb_* = \ml_1\, t_* \,,
$$
with $\ml_1=\ml_1(\xi,\xi_p^*,p)=\ml_{\min}(\Hb)$, the minimum eigenvalue of $\Hb$. Notice that $\ml_1>0$. $\ml_1$ depends on $\Mb_*$ which is unknown. Below we shall construct a lower bound $\underline{\ml_1}$ on $\ml_1$ and thus obtain a necessary condition for support points $\xb_*$ of $\xi_p^*$, in the form:
\begin{equation}\label{NC1}
    \xb_*\TT \Mb^{-(p+1)}\xb_* \geq  \underline{\ml_1}\, t_* \,.
\end{equation}

When $p=0$ ($D$-optimal design), we have $t=t_*=m$, and this necessary condition is simply
\begin{equation}\label{NC-p=0}
    \xb_*\TT \Mb^{-1}\xb_* \geq  \underline{\ml_1}\,m \ \ \ (p=0) \,;
\end{equation}
it corresponds to the case treated in \citet{HPa06_SPL}. When $p\neq 0$, $t_*$ is usually unknown and we shall use
\begin{eqnarray}
\xb_*\TT \Mb^{-(p+1)}\xb_* &\geq&  \underline{\ml_1}\, t(1+\me/t)^{-p} \ \ \  \mbox{ for } p>0 \,, \label{NC-p>0} \\
\xb_*\TT \Mb^{-(p+1)}\xb_* &\geq&  \underline{\ml_1}\, t  \ \ \ \hspace{1.8cm} \mbox{ for } -1<p<0 \,, \label{NC-p<0}
\end{eqnarray}
see (\ref{ineq-phi-phi_*-me}) and the definitions of $t, t_*, \phi_p, \phi_p^*$. Next section is devoted to the construction of the lower bound $\underline{\ml_1}$, using the inequalities (\ref{ineq1-t+me}) and (\ref{ineq2-t_*}).

%--------------------------------------------------------------------------------------
\subsection{Construction of the lower bound $\underline{\ml_1}$}

The inequality (\ref{ineq1-t+me}) can be rewritten as $\tr(\Hb^{1/(p+1)}\Mb^{-p}) \leq t+\me$ and (\ref{ineq2-t_*}) can be rewritten as $\tr(\Hb^{-1}\Mb^{-p}) \leq t_*$. Consider the spectral decomposition $\Hb=\Sb\mLb\Sb\TT$, with $\Sb\Sb\TT=\Sb\TT\Sb=\Ib_m$ and $\mLb$ the diagonal matrix whose diagonal elements are the eigenvalues $\ml_i$ of $\Hb$ sorted by increasing values. Denote $\Bb=\Sb\TT\Mb^{-p}\Sb$ and $b_i=\{\Bb\}_{ii}$ its diagonal elements, $i=1,\ldots,m$. $\Bb$ has the same set of eigenvalues as $\Mb^{-p}$ and
\begin{equation}\label{b_i}
    0< \underline{b_1} = \ml_{\min}(\Mb^{-p}) \leq b_i \leq \ml_{\max} (\Mb^{-p})\,, \ i=1,\ldots,m \,,
\end{equation}
as a consequence of Poincar\'e's separation Theorem, see, e.g., \citet[p.~211]{MagnusN99}.
We then obtain that (\ref{ineq1-t+me}) and (\ref{ineq2-t_*}) are respectively equivalent to
\begin{eqnarray}
\sum_{i=1}^m b_i\,\ml_i^{1/(p+1)} &\leq& t+\me \,, \nonumber \\
\sum_{i=1}^m b_i/\ml_i &\leq& t_* \,. \label{ineq2-ml}
\end{eqnarray}

\begin{rmk}\label{R:R1}
Inequality \eqref{ineq2-ml} implies that $\underline{\ml_1} \geq b_1/t_* \geq \underline{b_1}/t_*$. When plugged in \eqref{NC1}, it gives $\xb_*\TT \Mb^{-(p+1)}\xb_* \geq \underline{b_1}$. Although this bound is rather loose for $m\geq 2$, it cannot be improved when $m=1$. Indeed, $m=1$ implies $\underline{b_1}=b_1=t$ and the inequality $\xb_*\TT \Mb^{-(p+1)}\xb_* \geq  t$ is the tightest we can obtain, see Th.~\ref{Thm:EQ}-($ii$). In the following we shall suppose that $m\geq 2$.
\end{rmk}

Denote $\mo_i=\ml_i^{1/(p+1)}$ for $i=1,\ldots,m \geq 2$. The determination of $\underline{\ml_1}$ amounts to the solution of the following optimization problem: minimize $\mo_1$ with respect to $\mo=(\mo_1,\ldots,\mo_m)\TT$ under the constraints $0\leq \mo_1 \leq \mo_2 \leq \cdots \leq \mo_m$, $\sum_{i=1}^m b_i\,\mo_i \leq t+\me$ and $\sum_{i=1}^m b_i/\mo_i^{p+1} \leq t_*$. This is a convex problem, with Lagrangian
$$
L(\mob,\mu_1,\mu_2)=\mo_1+\mu_1\left( \sum_{i=1}^m b_i\,\mo_i - t -\me \right) + \mu_2 \left( \sum_{i=1}^m b_i/\mo_i^{p+1} - t_* \right)\,.
$$
Its stationarity with respect to $\mo$ indicates that the optimal solution satisfies $\mo_i=\mo_2$ for all $i\geq 2$. Since $\sum_{i=1}^m b_i=\tr(\Mb^{-p})=t$, from the Kuhn-Tucker conditions we obtain
\begin{eqnarray*}
b_1\,\mo_1 + (t-b_1) \mo_2 &=& t+\me \,,\\
b_1/\mo_1^{p+1} + (t-b_1)/ \mo_2^{p+1} &=& t_* \,,
\end{eqnarray*}
or equivalently
\begin{eqnarray}
\ma\,\mo_1 + (1-\ma) \mo_2 &=& 1+\beta \,, \label{bound1a} \\
\ma/\mo_1^{p+1} + (1-\ma)/ \mo_2^{p+1} &=& \mg_* \,, \label{bound1}
\end{eqnarray}
where $\ma=b_1/t$, $\beta=\me/t\geq 0$ and $\mg_*=t^*/t$.

When $p=0$ ($D$-optimal design), then $\ma=1/m$, $\mg_*=1$ and (\ref{bound1a}), (\ref{bound1}) can be directly solved for $\mo_1, \mo_2$, yielding $\underline{\ml_1}=\mo_1$ to be used in (\ref{NC-p=0}), see \citet{HPa06_SPL}. However, when $p\neq 0$, $\ma$ depends on $\Mb_*$ and $\mg_*$ depends on $t_*$ and are thus usually unknown. We must then determine the lowest value of $\mo_1\leq \mo_2$ satisfying (\ref{bound1a}), (\ref{bound1}) given the information available on $\ma$ and $\mg_*$; that is, respectively, (\ref{b_i}) which gives $1>\ma \geq \underline{b_1}/t = \ml_{\min}(\Mb^{-p})/\tr(\Mb^{-p})$, and (\ref{ineq-phi-phi_*-me}) which implies that $\mg_*$ satisfies
\begin{eqnarray}
\mg_*\in[(1+\beta)^{-p},1] && \mbox{ if } p\geq 0 \,, \label{mg-p>0} \\
\mg_*\in[1,(1+\beta)^{-p}] && \mbox{ if } p\leq 0 \,. \label{mg-p<0}
\end{eqnarray}

The solution to this problem is given in appendix and yields the main result of the paper.

\begin{thm}\label{Th:main}
For any $p\in(-1,\infty)$ and any design $\xi\in\Xi^+$, any point $\xb_*\in\SX$ such that
\begin{equation}\label{remove-condition}
    \xb_*\TT \Mb^{-(p+1)}(\xi) \xb_* < C(\xi,p)= \mo_1^{p+1} \, B(t,\me)
\end{equation}
cannot be support point of a $\phi_p$-optimal design measure $\xi_p^*$, where we denoted $t=\tr[\Mb^{-p}(\xi)]$, $\me=\max_{\xb\in\SX} \xb\TT \Mb^{-(p+1)}(\xi)\xb -t$, $B(t,\me)=t\,\min\{1,(1+\me/t)^{-p}\}$, and where $\mo_1$ is the unique solution for $\mt$ in the interval
$((\ma/\mg)^{1/(p+1)},(1/\mg)^{1/(p+1)}]$ of the equation
\begin{equation}\label{F}
    F(\mt; \ma,\me, t,\mg,p) = \frac{\ma}{\mt^{p+1}} + \frac{(1-\ma)^{p+2}}{(1+\me/t-\ma\mt)^{p+1}} - \mg =0
\end{equation}
with $\ma=\ml_{\min}[\Mb^{-p}(\xi)]/\tr[\Mb^{-p}(\xi)]$ and $\mg=\max\{1,(1+\me/t)^{-p}\}$.

In the special case when $t_*=\tr[\Mb^{-p}(\xi_p^*)]$ is known (thus, in particular if $p=0$), one can take $B(t,\me)=t_*$ and $\mg=\mg_*=t^*/t$ in \eqref{remove-condition} and \eqref{F}.
\end{thm}

Denote $\delta=\max_{\xb\in\SX} F_{\phi_p}(\xi,\xb)$. The theorem indicates that any support point $\xb_*$ of a $\phi_p$-optimal design measure satisfies the inequality
$
F_{\phi_p}(\xi,\xb_*) \geq h_p[\Mb(\xi),\delta]\,,
$
where $h_p[\Mb(\xi),\delta]=\phi_p(\xi)\,[\mo_1^{p+1} \, B(t,\me)/t -1]$ with $\me=\delta t/\phi_p(\xi)$, see \eqref{F(xi,x)}. Notice that $\mo_1^{p+1} \, B(t,\me)\leq t$, so that $h_p[\Mb(\xi),\delta]\leq 0$ and all points $\xb$ such that $F_{\phi_p}(\xi,\xb) \geq 0$ are potential support points of $\xi_p^*$. When $\delta$ tends to zero, then $\me \ra 0$ and $h_p[\Mb(\xi),\delta] \ra 0$, see the proof of Th.~\ref{Th:main}, in accordance with the last statement of the Equivalence Theorem.

\begin{rmk} \mbox{}\\
\vspace{-1cm}
\begin{enumerate}
  \item When $p$ is integer, $F(\mt; \ma,\me,t,\mg,p)=0$ is a polynomial equation in $\mt$ of degree $2(p+1)$.

  \item Suppose $p>0$ with $t_*$ unknown and $\me\ra\infty$; then, $B(t,\me)\ra 0$, so that $C(\xi,p)\ra 0$ and the condition \eqref{remove-condition} brings no information on the support of $\xi_p^*$. The same is true when $p<0$ with $t_*$ unknown and $\me\ra\infty$: $\mg\ra\infty$, so that $\mo_1\ra 0$ and again $C(\xi,p)\ra 0$. Suppose now that $t_*$ is known. Then, $C(\xi,p)=t_*\, \mo_1^{p+1} \in (\ml_{\min}[\Mb^{-p}(\xi)],\, \tr[\Mb^{-p}(\xi)]]$ and $\mo_1^{p+1} \ra \ma/\mg_*=\ml_{\min}[\Mb^{-p}(\xi)]/t_*$ as $\me\ra\infty$, see \eqref{F}, so that $C(\xi,p)\ra \underline{b_1}= \ml_{\min}[\Mb^{-p}(\xi)]$ and we recover the same bound as in Remark~\ref{R:R1}.

  \item Using a construction similar to that in \citet[Th.~3]{HPa06_SPL}, one can show that the bound \eqref{remove-condition} with  $B(t,\me)=t_*$ and $\mg=t^*/t$ gives the tightest necessary condition for support points: for any $m\geq 2$, any $\me,\, \me'>0$, one can exhibit an example with a design space $\SX$, a design measure $\xi$ such that $\max_{\xb\in\SX} \{\xb\TT \Mb^{-(p+1)}\xb\} - t=\me$, and an optimal design $\xi_p^*$ with support point $\xb_*$ such that $\xb_*\TT \Mb^{-(p+1)}\xb_* < \mo_1^{p+1}  \,t_*+ \me'$ (with $\Mb$ and $\Mb_*$ diagonal and $\Hb$ having eigenvalues $\ml_1 <\ml_2 = \cdots = \ml_m$).

\end{enumerate}
%1.
%\noindent 2.
\end{rmk}

%______________________________________________________________________________________
\section{Examples}\label{S:Example}

\paragraph{Example 1.}

\begin{figure}
\begin{center}
\includegraphics[width=.45\linewidth]{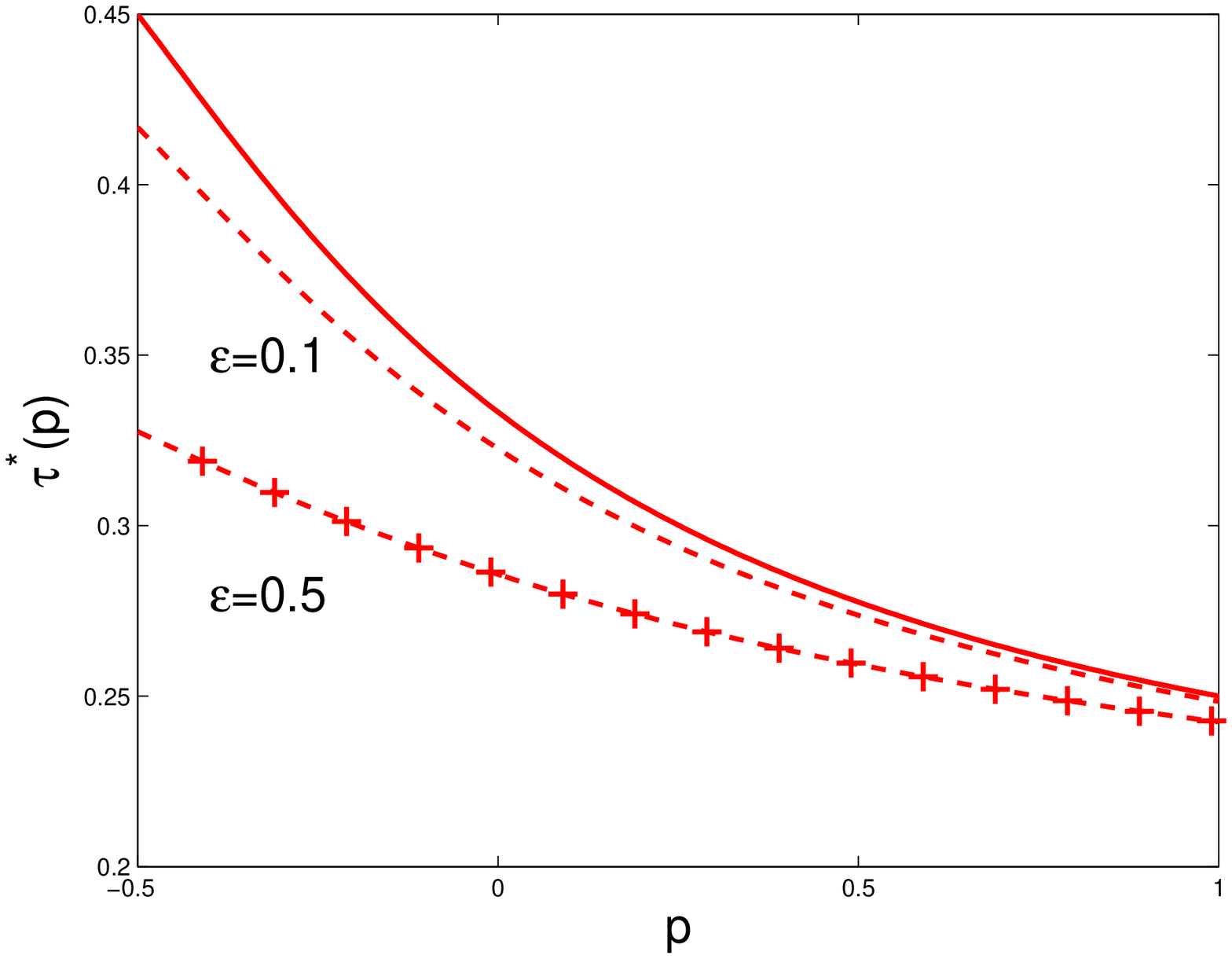}
\includegraphics[width=.45\linewidth]{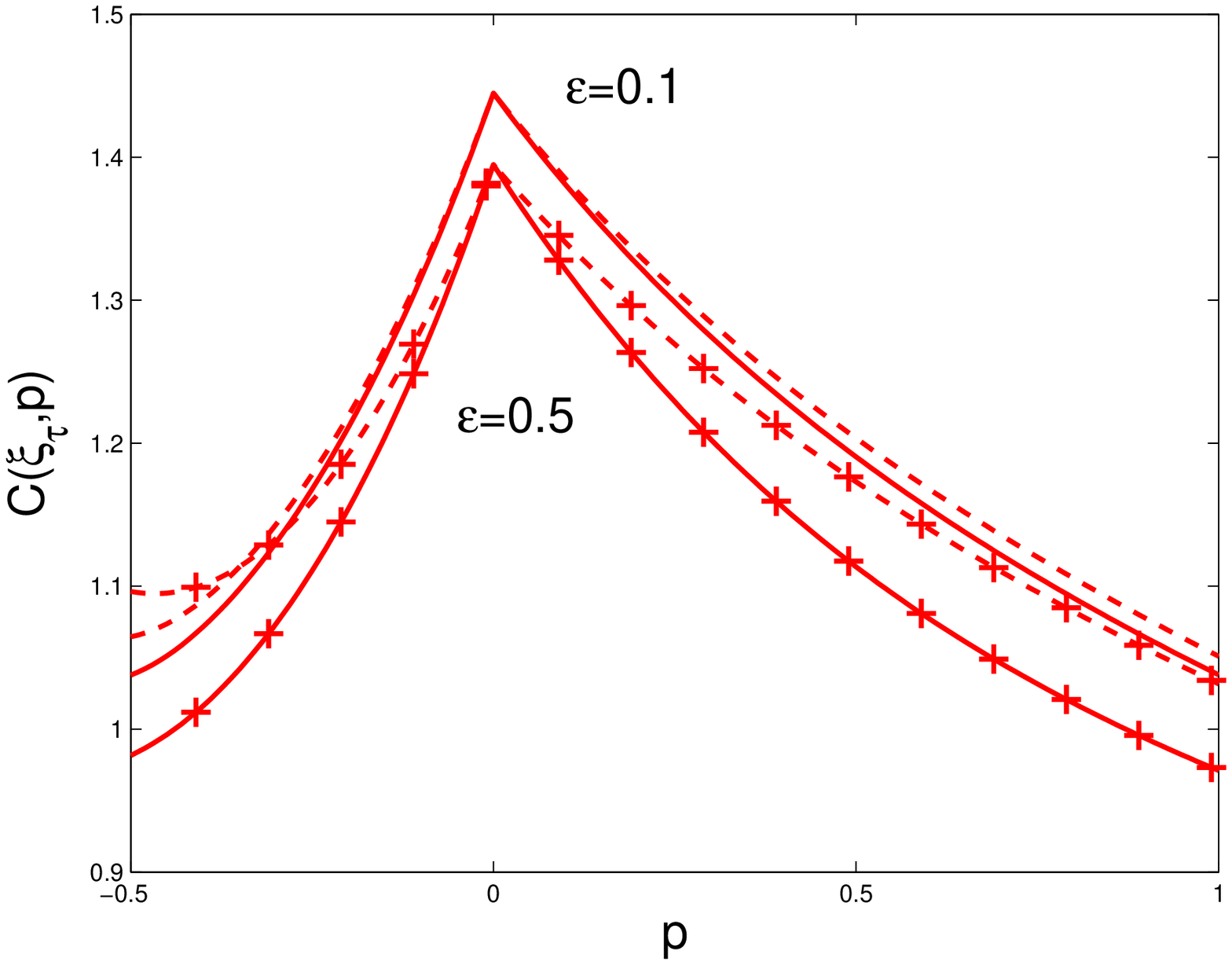}
\end{center}
\caption{\footnotesize Left: $\tau^*(p)$ such that $\xi_{\tau^*(p)}=\xi_p^*$ is $\phi_p$-optimal for $p$ (\emph{solid line}) and
$\tau(p,\me)$ such that $\max_{s\in[-1,1]} \xb\TT(s) \Mb^{-(p+1)}(\xi_{\tau(p,\me)}) \xb(s) = \tr[\Mb^{-p}(\xi_{\tau(p,\me)})] + \me$ ($\me=0.1$ and $\me=0.5$, \emph{dashed lines}). Right: bound $C(\xi_{\tau(p,\me)},p)$ in (\ref{remove-condition}) for the two cases $t_*$ unknown (\emph{solid lines}) and $t_*$ known (\emph{dashed lines}) for $\me=0.1$ and $\me=0.5$. The curves for $\me=0.5$ are marked with crosses.}
\label{F:tau_bounds}
\end{figure}

Consider the linear regression model with $\xb=\xb(s)=(1,\ s,\ s^2)\TT$, $s\in[-1,1]$ ($m=3$). For any $p\in(-1,\infty)$, the $\phi_p$-optimal design on $[-1,1]$ is unique and is supported at the three points $\{-1,0,1\}$. For symmetry reasons, it corresponds to
$$
\xi_\tau = \tau\, \delta_{-1} + (1-2\tau)\, \delta_0 + \tau\,\delta_1
$$
for some particular $\tau^*=\tau^*(p)$, with $\tau^*(-1/2)=0.45$, $\tau^*(0)=1/3$ ($D$-optimal design), $\tau^*(1)=1/4$ ($A$-optimal design) and, in the limit $p\ra\infty$, $\tau^*(\infty)=0.2$ ($E$-optimal design), see Fig.~\ref{F:tau_bounds}-left for a plot of $\tau^*(p)$ for $p\in[-1/2,1]$. Here, $\delta_s$ denotes the Dirac delta measure at $s$.

\begin{figure}
\begin{center}
\includegraphics[width=.5\linewidth]{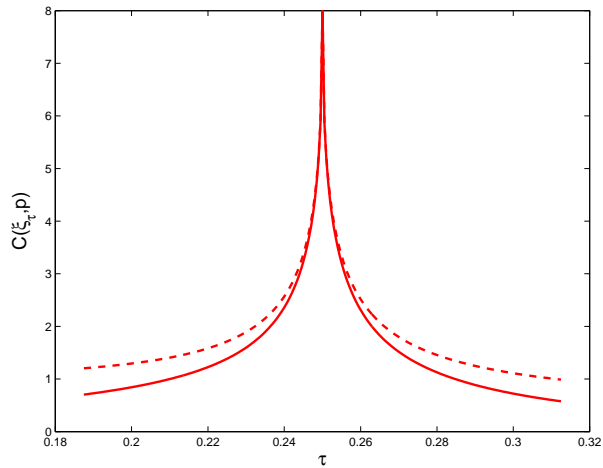}
\end{center}
\caption{\footnotesize Value of $C(\xi_\tau,p)$ for different designs $\xi_\tau$, $\tau\in[3/16,\, 5/16]$ ($p=1$, $t_*$ unknown in \emph{solid line}, $t_*$ known in \emph{dashed line}).}
\label{F:L1_tau}
\end{figure}

To illustrate the impact of not knowing $t_*$ on the construction of $\omega_1^{p+1}$ through the solution of (\ref{F}), we take $p=1$ and compute the bound $C(\xi_\tau,p)$ for the cases $\mg=1$ ($t_*$ unknown) and $\mg=t^*/t$ ($t_*$ known) for different designs $\xi_\tau$, $\tau\in[\tau^*(1)-1/16,\, \tau^*(1)+1/16]$. Figure~\ref{F:L1_tau} shows that, close to the optimum $\tau^*(1)=1/4$, the value obtained for $t_*$ unknown (solid line) is not much worse, i.e., smaller, than the value for $t_*$ known (dashed line). Note that considering different designs $\xi_\tau$ with $\tau\neq \tau^*(p)$ is equivalent to considering different $\me$ given by (\ref{epsilon}), with $\me$ being approximately linear in $|\tau-\tau^*(p)|$ for the range of values of $\tau$ considered.

The marginal deterioration of the bound (\ref{remove-condition}) due to the ignorance of $t_*$ when $\me$ is small enough is further illustrated by Fig.~\ref{F:tau_bounds}. Here, we set $\me$ at some fixed value (the values $\me=0.1$ and $\me=0.5$ are considered), and for values of $p$ in the range $[-1/2,1]$ we compute $\tau(p,\me)$ such that $\max_{s\in[-1,1]} \xb\TT(s) \Mb^{-(p+1)}(\xi_\tau) \xb(s) = \tr[\Mb^{-p}(\xi_\tau)] + \me$. The values of $\tau^*(p)$ and $\tau(p,\me)$, $\me=0.1$, $0.5$, are shown in Fig.~\ref{F:tau_bounds}-left, in solid and dashed lines respectively. Then, for each $p$ and associated design $\xi_{\tau(p,\me)}$ we compute the bound $C(\xi_{\tau(p,\me)},p)$ of (\ref{remove-condition}) in the two situations $t_*$ unknown and $t_*$ known; see the plots in Fig.~\ref{F:tau_bounds}-right. Note that the bound for $t_*$ unknown (solid line) remains near the bound for $t_*$ known (dashed line) when $\me=0.1$; the situation deteriorates for larger $\me$ (curves with crosses) but the two bounds get close as $p$ approaches 0 and exactly coincide at $p=0$ (since then $t=t_*=m$).

\paragraph{Example 2.}

Take now the complete product-type interaction model $\xb(\sb)=\xb(s_1) \otimes \xb(s_2)$, $\sb=(s_1,\ s_2)$,  with $\otimes$ denoting tensor product and $\xb(s_i)=(1,\ s_i,\ s_i^2)\TT$, $s_i\in[-1,1]$, for $i=1,2$ ($m=9$). The $D$-optimal (respectively $A$-optimal) design for this problem is the cross product of two $D$-optimal designs (resp.\ $A$-optimal designs) for one single factor, i.e., it corresponds to the cross product of two designs $\xi_\tau$ with $\tau=1/3$ (resp.\ $\tau=1/4$), see \citet[Chap.~4 and 5]{Schwabe96}. The optimal values of $\phi_p$, $p=0,1$, are $\phi_0^*=16^{1/3}/9 \simeq 0.2800$ and $\phi_1^*=9/64 \simeq 0.1406$.

\begin{figure}
\begin{center}
\includegraphics[width=.45\linewidth]{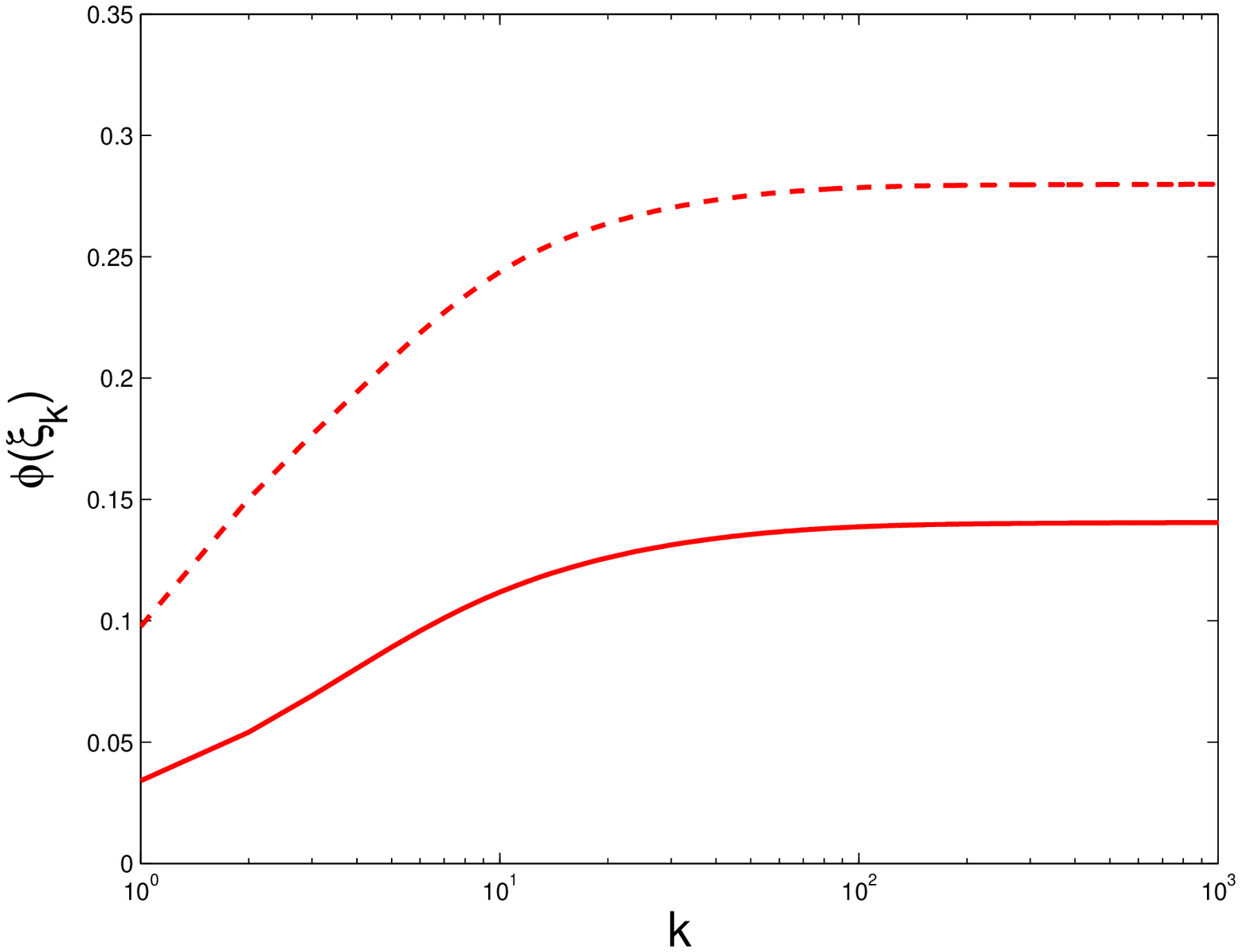}
\includegraphics[width=.45\linewidth]{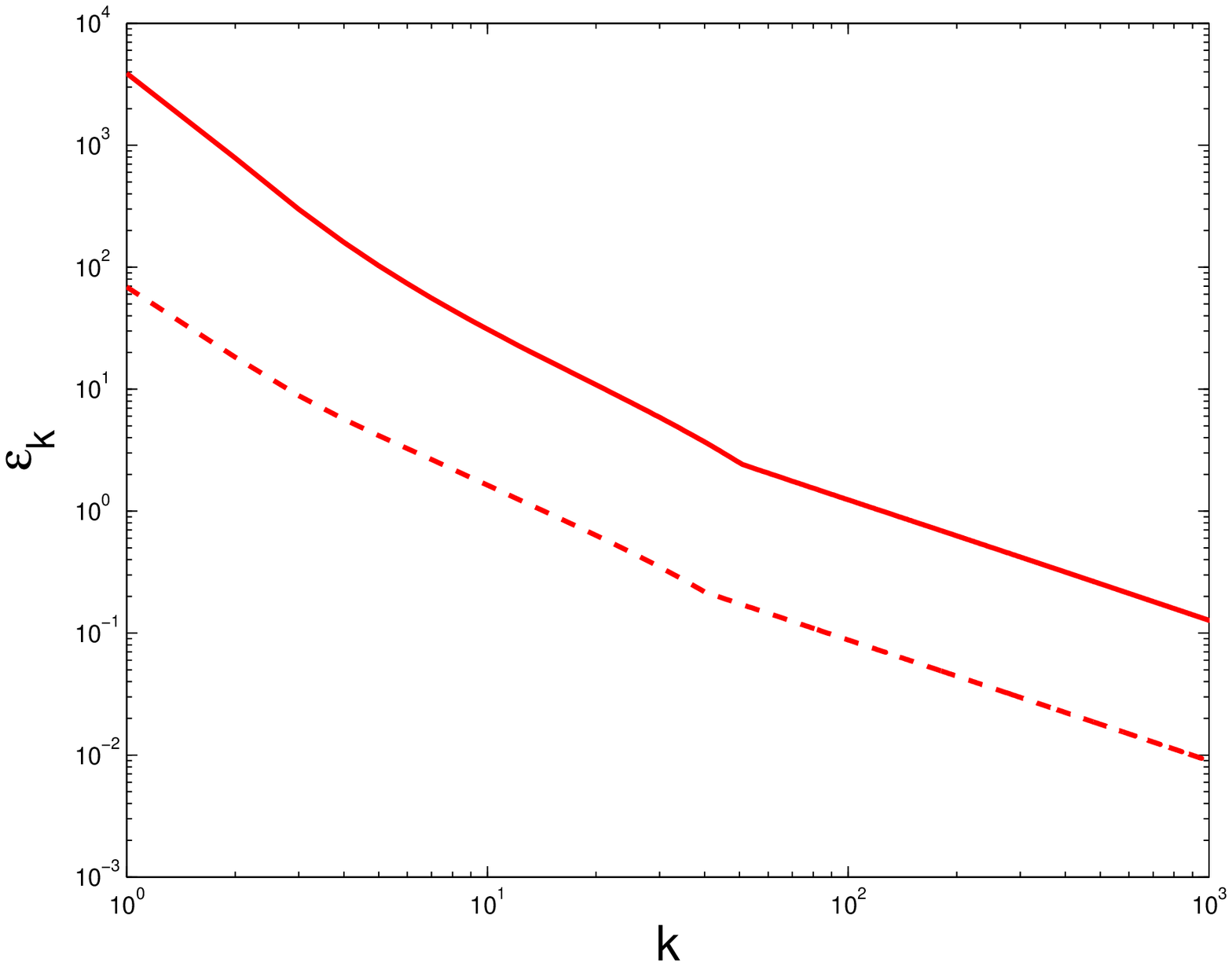}
\end{center}
\caption{\footnotesize $\phi(\xi_k)$ --- left --- and $\me_k=\me(\xi_k)$ given by (\ref{epsilon}) --- right --- as functions of $k$ for the recursion (\ref{iterations}); $D$-optimal design is in \emph{dashed line}, $A$-optimal design is in \emph{solid line}. }
\label{F:phik-epsik}
\end{figure}

We consider the iterative construction of optimal designs through the recursion
\begin{equation}\label{iterations}
    w_i^{k+1} = w_i^k \frac{ [\xb_i\TT \Mb^{-(p+1)}(\xi_k)\xb_i]^a }{\sum_{i=1}^{N_k} [\xb_i\TT \Mb^{-(p+1)}(\xi_k)\xb_i]^a } \,,
\end{equation}
where $\xi_k$, the design measure at iteration $k$, allocates mass $w_i^k$ at the point $\xb_i$ present in $\SX$ at iteration $k$, $i=1,\ldots,N_k$. The initial design space corresponds to a uniform grid for $\sb$, with $s_i$ varying from $-1$ to $1$ by steps of $0.01$ ($201$ values), $i=1,2$, which gives $N_0=40~401$. The initial design $\xi_0$ is the uniform measure on those $N_0$ points. We take $a=1$ for $D$-optimal design ($p=0$) and $a=1/2$ for $A$-optimal design ($p=1$), which ensures monotonic convergence to the optimum, see \citet{Titterington76} and \citet{Pazman86} for $D$-optimal design and \citet{Torsney83} for $A$-optimal design; see also Fig.~\ref{F:phik-epsik}-left. One may also refer to \citet{SilveyTT78} for a general class of multiplicative algorithms and to \citet{DettePZ2008} for an improved updating rule yielding accelerated convergence.
Due to the convergence of $\xi_k$ to the optimal design, $\me_k=\me(\xi_k)$ given by (\ref{epsilon}) is decreasing with $k$, see Fig.~\ref{F:phik-epsik}-right.

\begin{figure}
\begin{center}
\includegraphics[width=.5\linewidth]{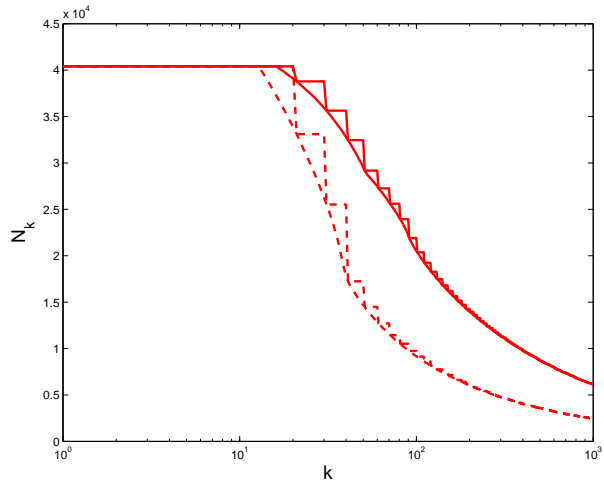}
\end{center}
\caption{\footnotesize $N_k$ as a function of $k$ when using (\ref{remove-condition}) to remove points from $\SX$; for $D$-optimal design (\emph{dashed line}) and $A$-optimal design (\emph{solid line}). The condition (\ref{remove-condition}) is used every iteration for the continuous curves, every 10th iteration for the staircase curves.}
\label{F:Nk}
\end{figure}

We use inequality (\ref{remove-condition}) to reduce the cardinality $N_k$ of $\SX$ when possible: any point that violates (\ref{remove-condition}) cannot be a support point of the optimal measure and is removed from $\SX$. Here we simply set its mass to zero and rescale the weights of remaining point so that they sum to one, but more sophisticated reallocation rules can be used, see \citet{HPa06_SPL}. $N_k$ thus decreases with $k$, rendering the iterations (\ref{iterations}) simpler and simpler as $k$ increases. Figure~\ref{F:Nk} shows the evolution of $N_k$ with $k$, both for $D$-optimal and $A$-optimal designs (in dashed and solid line respectively): cancelation of points is performed at every iteration for the continuous curves, every 10th iterations only for the staircase curves.

The decrease of $N_k$ is faster for $D$-optimal design than for $A$-optimal design, the bound $C(\xi,p)$ in (\ref{remove-condition}) being more pessimistic for the latter, see Fig.~\ref{F:tau_bounds}-right, and $\me$ being larger, see Fig.~\ref{F:phik-epsik}-right.
Note that the cancelation of points does not hamper the convergence of (\ref{iterations}) since (\ref{remove-condition}) is used a finite number of times only (obviously bounded by $N_0$) --- the heuristic rule used to reallocate weights of points that are removed may, however, impact monotonicity, although this is not the case in the present example, see Fig.~\ref{F:phik-epsik}-left. Also, the effect of cancelation on the behaviors of $\phi(\xi_k)$ and $\me_k=\me(\xi_k)$ (Fig.~\ref{F:phik-epsik}) is negligible: the acceleration of the algorithm is only due to the reduction of the cardinality $N_k$. Taking as reference $t_c=1$ the computing time for 1~000 iterations of the recursion (\ref{iterations}) for $D$-optimal design with cancelation of points at each iteration, we get $t_c=11.5$ for $D$-optimal design without cancelation, and $t_c=12.15$, $t_c=2.8$, for $A$-optimal design, respectively without and with cancelation at each iteration. Cancelation need not be checked at each iteration though, and the computing times become $t_c=0.87$ and $t_c=2.4$ for $D$- and $A$-optimal designs respectively when the condition (\ref{remove-condition}) is used each 10th iteration only (see the staircase curves on Fig.~\ref{F:Nk}). Clearly, a suitable adaptation of the frequency of cancelation of points
to the decrease of $N_k$ might provide further reductions in computing time.

%______________________________________________________________________________________
\section{Possible extensions and conclusions}\label{S:Extensions}

Multivariate regression and Bayesian optimal design involve information matrices that can be expressed as $\Mb(\xi)=\int_\SX \SM(\xb)\,\xi(\dd\xb)$ with $\SM(\xb)\in\mathbb{M}$ having rank larger than one (we suppose that $\SM(\cdot)$ is measurable and that $\{\SM(\xb), \xb\in\SX\}$ forms a compact subset of $\mathbb{M}$). The results presented here can easily be extended to that situation, following the same lines as in \citet{HarmanT2009} where the case $p=0$ is considered.

The $E$-optimality criterion $\phi_E(\xi)=\Phi_E[\Mb(\xi)]=\ml_{\min}[\Mb(\xi)]$ is not differentiable in general, but $\Phi_E(\cdot)$ is differentiable at $\Mb$ when $\ml_{\min}(\Mb)$ has multiplicity one, with gradient $\nabla \phi_E(\Mb)=\vb\vb\TT$ where $\vb$ denotes the eigenvector of unit length (unique up to a sign change) associated with $\ml_{\min}(\Mb)$. Although $\phi_E(\xi)$ corresponds to the limit of $\phi_p^+(\xi)$ as $p$ tends to infinity, the results of Sect.~\ref{S:Phi_p} do not extend to this limiting situation, even in the differentiable case; $E$-optimality thus requires a special treatment and will be considered elsewhere.

The determination of a $D$-optimal design can be used for maximum-likelihood estimation in mixture models, see, e.g., \citet{Lindsay83} and \citet{Mallet86}, and for the construction of the minimum-volume ellipsoid containing a compact set, see, e.g., \citet{Sibson72}, \citet{KhachiyanT93} and \citet{Khachiyan96}. More generally, for any $q\in(-1,\infty)$ the determination of the ellipsoid $\SE(\Ab)=\{\zb\in\mathbb{R}^m: \zb\TT \Ab\zb \leq 1\}$, $\Ab\in\mathbb{M}$, containing the $k$ points $\xb_1,\ldots,\xb_k$ of $\mathbb{R}^m$ and such that $\phi_q(\Ab)$ is maximum is equivalent to the determination of a $\phi_p$-optimal design on $\SX=\{\xb_1,\ldots,\xb_k\}$ with $p=-q/(1+q)\in(-1,\infty)$, and the optimal matrix $\Ab_*$ equals
$\Mb_*^{-(p+1)}/t_*$; see \citet[Chap.~6]{Pukelsheim93}. The delimitation of the support points of a $\phi_p$-optimal design can therefore also be used to accelerate the algorithmic construction of ``$\phi_q$-optimal ellipsoids'' containing compact sets. Note that, as illustrated in \citep{Pa03, HPa06_SPL}, a more substantial acceleration than in the optimal design example of Sect.~\ref{S:Example} can be expected.

In Sect.~\ref{S:Example}, we considered the suppression of points that cannot be support points of an optimal design in a multiplicative algorithm. When $\SX$ is not finite, or is finite but very large, it is advisable to use a vertex-direction or a vertex-exchange algorithm, see, e.g., \citet{Fedorov72}, \citet{Wu78} and \citet{Bohning86}. This requires the determination at each iteration, say iteration $k$, of a point $\hat\xb_k$ of $\SX$ that maximizes $F_{\phi_p}(\xi_k,\xb)$ given by (\ref{F(xi,x)}), at least approximately.
Condition (\ref{remove-condition}) of Theorem~\ref{Th:main} can then be used to restrict the search for a suitable $\hat\xb_k$ in a domain that shrinks as $k$ increases. Further developments are required to construct algorithms making an efficient use of (\ref{remove-condition}) for the inclusion of new support points.

%______________________________________________________________________________________
\section*{Appendix}

\noindent{\em Proof of Th~\ref{Th:main}.} %\mbox{}\\
The proof is in three parts. In ($i$) we show that for given $\ma$ and $\mg_*$ the equations (\ref{bound1a}), (\ref{bound1}) with $\mo_1\leq \mo_2$ have a unique solution $\mo_1^*(\ma,\mg_*)$ for $\mo_1$, with $\mo_1^*(\ma,\mg_*)\in((\ma/\mg_*)^{1/(p+1)},(1/\mg_*)^{1/(p+1)}]$. Then in ($ii$) we show that this solution is non-decreasing in $\ma$, so that the required lowest bound is obtained for $\ma=\underline{b_1}/t$, see (\ref{b_i}). Finally, in ($iii$) we consider the case when $t_*$ is unknown.

($i$) Expressing $\mo_2$ as a function of $\mo_1$ using (\ref{bound1a}), we obtain $\mo_2=f_1(\mo_1)=(1+\beta-\ma\mo_1)/(1-\ma)$, i.e., a decreasing linear function of $\mo_1$ with slope $-\ma/(1-\ma)$ and such that $f_1[(1+\beta)/\ma]=0$. Doing the same with (\ref{bound1}), we obtain $\mo_2=f_2(\mo_1)$ with $f_2(\cdot)$ decreasing and concave for $\mo_1 \in((\ma/\mg_*)^{1/(p+1)},\infty)$, $f_2(\mt)$ tending to infinity when $\mt$ approaches $(\ma/\mg_*)^{1/(p+1)}$ from above and $\lim_{\mt\ra\infty} f_2(\mt)=1/\ma-1$. Note that (\ref{mg-p>0}), (\ref{mg-p<0}) imply that $(\ma/\mg_*)^{1/(p+1)}<(1/\mg_*)^{1/(p+1)}<(1+\beta)/\ma$. Therefore, $f_2(\mt)>f_1(\mt)$ for $\mt$ close enough to $(\ma/\mg_*)^{1/(p+1)}$ or large enough.

Denote $f_2'(\mt)=\dd f_2(\mt)/\dd \mt$ and consider $\mt_*=(1/\mg_*)^{1/(p+1)}$. Direct calculations indicate that $f_2(\mt_*)=\mt_*$, $f_2'(\mt_*)=-\ma/(1-\ma)$ with, moreover, $f_1(\mt_*) > f_2(\mt_*)$ when $\beta>0$, i.e., when $\me>0$, due to (\ref{mg-p>0}) and (\ref{mg-p<0}). Two solutions $\mo_{1,a}^*, \mo_{1,b}^*$ thus exist for (\ref{bound1a}), (\ref{bound1}), with $\mo_{1,a}^* < \mt_* < \mo_{1,b}^*$. Only $\mo_{1,a}^*$ is such that the associated $\mo_{2,a}^*$ satisfies
$\mo_{2,a}^*> \mo_{1,a}^*$. When $\me=0$, then $f_1(\mt_*) = f_2(\mt_*) =\mt_*$ and the two solutions $\mo_{1,a}^*, \mo_{1,b}^*$ are confounded and equal $\mt_*$ (and also coincide with $\mo_{2,a}^*$ and $\mo_{2,b}^*$).
The equations (\ref{bound1a}) and (\ref{bound1}) with $\mo_1\leq \mo_2$ thus always have a unique solution $\mo_1^*(\ma,\mg_*)$ and this solution belongs to the interval $((\ma/\mg_*)^{1/(p+1)},\mt_*]$.

($ii$) Applying the implicit function theorem to (\ref{bound1a}), (\ref{bound1}) we obtain that the solution $\mo_1^*(\ma,\mg_*)$ satisfies
\begin{eqnarray*}
\frac{\mp \mo_1^*(\ma,\mg_*)}{\mp\ma} &=& \frac{(p+1) (\mo_1^*)^{p+2}(\mo_1^*-\mo_2^*)+\mo_1^*\mo_2^* [(\mo_2^*)^{p+1}-(\mo_1^*)^{p+1}]}{\ma(p+1)[(\mo_2^*)^{p+2}-(\mo_1^*)^{p+2}]} \\
&=& \frac{\mo_1^*}{\ma(p+1)(z^{p+2}-1)}\ [(p+1)(1-z)+z(z^{p+1}-1)] \,,
\end{eqnarray*}
where $z=\mo_2^*/\mo_1^* \geq 1$. Denote $f(z)=(p+1)(1-z)+z(z^{p+1}-1)$, its derivative is $\dd f(z)/\dd z=(p+2) (z^{p+1}-1)$ so that $f(z) \geq f(1)=0$.
Since (\ref{b_i}) gives $\ma \geq \underline{b_1}/t$, one has $\mo_1^*(\ma,\mg_*) \geq \mo_1^*(\underline{b_1}/t,\mg_*)$. The substitution of $[\mo_1^*(\underline{b_1}/t,\mg_*)]^{p+1}$ for $\underline{\ml_1}$ in (\ref{NC1}) concludes the proof for the case when $t_*$ is known.

($iii$) When $t_*$ is unknown, an upper bound can be substituted for $t_*$ in (\ref{ineq2-ml}). Using (\ref{mg-p>0}), (\ref{mg-p<0}), this amounts at replacing $\mg_*$ by the upper bound $\mg=\max\{1,(1+\me/t)^{-p}\}$. The necessary conditions (\ref{NC-p>0}), (\ref{NC-p<0}) with $\underline{\ml_1}=[\mo_1^*(\underline{b_1}/t,\mg)]^{p+1}$ then give (\ref{remove-condition}).

%______________________________________________________________________________________
\section*{Acknowledgments}
The author thanks the two referees and the associate editor for their constructive comments on an earlier version of the paper.

%______________________________________________________________

%\bibliographystyle{plain}
%\bibliographystyle{unsrt}
\bibliographystyle{elsart-harv}
%\bibliography{xampl,test}

\end{document}